\documentclass[12pt,draft,twoside]{article} 
\usepackage{amsmath,amssymb,latexsym,theorem,bbm,shapepar,natbib}
\def\Sign{\mathop{\hbox{\rm sign}}}
    
\newcommand{\Var}{\mathrm{Var}}
\newcommand{\Cov}{\mathrm{Cov}}

\newcommand{\bone}{\mathbbm{1}}
\newcommand{\stoch}{\stackrel{\scriptstyle \mathsf P}{\longrightarrow}}
\newcommand{\distr}{\stackrel{\scriptstyle \mathcal D}{\longrightarrow}}
\newcommand{\qmean}{\stackrel{\scriptstyle {\mathsf L}_2}{\longrightarrow}}
\newcommand{\distre}{\stackrel{\scriptstyle \mathcal D}{=}}

\newcommand{\proofend}{\hfill$\square$}

\numberwithin{equation}{section}

\theorembodyfont{\em}
\newtheorem{Lem}{Lemma}[section]
\newtheorem{Thm}[Lem]{Theorem}
\newtheorem{Pro}[Lem]{Proposition}

\theorembodyfont{\rm}

\title{Testing stability in a spatial unilateral autoregressive model}

\author{\sc S\'andor Baran$^1$, Gyula Pap$^2$ and Kinga Sikolya$^1$ \\
$^1$Faculty of Informatics, University of Debrecen\\
         Kassai \'ut 26, H--4028 Debrecen, Hungary \\ [5mm]
$^2$Bolyai Institute, University of
             Szeged \\ Aradi v\'ertan\'uk tere 1, H-6720 Szeged, 
             Hungary 
}

\date{}

\begin{document}

\maketitle

\begin{abstract}
Least squares estimator of the stability parameter
 \ $\varrho := |\alpha| + |\beta|$ \ for a spatial unilateral autoregressive
 process \ $X_{k,\ell} = \alpha X_{k-1, \ell} + \beta X_{k, \ell-1} +
 \varepsilon_{k, \ell}$ \ is investigated.
Asymptotic normality with a scaling factor \ $n^{5/4}$ \ is shown in the
 unstable case, i.e., when \ $\varrho = 1$, \ in contrast to the AR(p) model
 \ $X_k = \alpha_1 X_{k-1} + \cdots + \alpha_p X_{k-p} +
 \varepsilon_k$, \ where the  least squares estimator of the stability parameter
 \ $\varrho := \alpha_1 + \cdots + \alpha_p$ \ is not asymptotically normal in
 the unstable, i.e., in the unit root case.



\end{abstract}

\section{Introduction}
  \label{sec:sec1}
Consider a spatial unilateral autoregressive process
\ $\{ X_{k,\ell} : k, \ell \in {\mathbb Z}, \, k + \ell \geq0 \}$ \ defined by
 \begin{equation}
   \label{model}
   X_{k,\ell} = \begin{cases}
               \alpha X_{k-1, \ell} + \beta X_{k, \ell-1} +
               \varepsilon_{k, \ell} , 
                & \text{for \ $k + \ell \geq 1$,}\\
               0,  & \text{for \ $k + \ell = 0$,}
              \end{cases}
 \end{equation}
 where \ $\{ \varepsilon_{k, \ell} : k, \ell \in {\mathbb Z}, \, k +
 \ell \geq 1 \}$ \ are
 independent random variables with \ ${\mathsf
   E}(\varepsilon_{k,\ell}) = 0$ \ and \ $\Var(\varepsilon_{k,\ell}) =
 1$. \ This model is stable in case of \ $|\alpha| + |\beta| < 1$ \  and
unstable if \ $|\alpha| + |\beta| = 1$ \ (see
\citet{whittle,besag,br2}), hence \ $\varrho := |\alpha|+ |\beta |$ \
can be considered as a {\em stability parameter\/}.   

For a set \ $H \subset \{(k, \ell) \in {\mathbb Z}^2 : k + \ell \geq
1\}$, \ the least 
 squares estimator (LSE) \ $(\widehat\alpha_H, \widehat\beta_H)$ \ of
 the coefficients \
 $(\alpha, \beta)$ \ based 
 on the observations \ $\{X_{k,\ell} : (k, \ell) \in H\}$ \ can be obtained by 
 minimizing the sum of squares
 \begin{equation*}
   \sum_{(k, \ell) \in H}
    \big( X_{k, \ell} - \alpha X_{k-1, \ell} - \beta X_{k, \ell-1} \big)^2
 \end{equation*}
 with respect to \ $\alpha$ \ and \ $\beta$, \ and it has the form
 \begin{equation*}
   \begin{bmatrix} \widehat\alpha^*_H \\ \widehat\beta^*_H \end{bmatrix}
   \!=\! (A_H^*)^{-1}b_H^*, \quad \text{where} \quad
A_H^*:=\!\!\!\sum_{(k, \ell) \in H}
            \begin{bmatrix}
             X_{k-1, \ell}    \\
             X_{k, \ell-1}
            \end{bmatrix} 
             \begin{bmatrix}
             X_{k-1, \ell} \\
             X_{k,\ell-1}
            \end{bmatrix}^{\top}\!\!\!, \quad 
     b_H^*:=\!\!\!\sum_{(k, \ell) \in H}\!\!\! X_{k, \ell}
      \begin{bmatrix}
       X_{k-1, \ell} \\
       X_{k, \ell-1}
      \end{bmatrix}.
 \end{equation*}

Model \eqref{model} has been investigated in details by several
authors. \citet{paul} determined the exact asymptotic behaviour of the
variances of the process, while \citet{bpz2} proved the asymptotic
normality of the LSE of the coefficients \ $(\alpha,\beta)$ \ both in
stable and unstable cases.

The limiting behavior of the LSE of the stability parameter \ $\varrho$ \ has
 not been treated yet, but such a stability parameter is well investigated
 in case of unstable AR($p$) processes, see the unit root tests, e.g., in
 \citet[Section 17, Table 17.3, Case 1]{hamilton}.
Namely, for the simplicity, in case of an AR(1) process
 \ $Y_k = \varrho Y_{k-1} +\zeta_k$, \ $k\in{\mathbb N}$, \ with \
 $Y_0:=0$ \ and an 
 i.i.d.\ sequence \ $\{\zeta_k : k\in{\mathbb N}\}$ \ having mean \
 $0$ \ and positive 
 variance, the LSE of the parameter \ $\varrho$ \ based on a sample
 \ $\{Y_1,\ldots,Y_n\}$ \ takes the form
\begin{equation*}
   \widehat\varrho_n = \frac{\sum_{k=1}^n Y_{k-1}Y_k}{\sum_{k=1}^n
     Y_k^2},\qquad n\in{\mathbb N},
 \end{equation*}
 see, e.g., \citet[17.4.2]{hamilton}, and, by \citet[17.4.7]{hamilton}, in the
 unstable case, i.e., when \ $\varrho=1$,
 \begin{equation*}
   n(\widehat\varrho_n-1) \distr \frac{\int_0^1 {\mathcal W}_t{\mathrm d}
     {\mathcal W}_t}{\int_0^1 {\mathcal W}_t^2{\mathrm d} t} 
   \qquad \text{as \ $n\to\infty$,}
 \end{equation*}
 where \ $\big({\mathcal W}_t\big)_{t\geq 0}$ \ is a standard Wiener process.
Here \ $n(\widehat\varrho_n-1)$ \ is called Dickey-Fuller unit root statistics.
It turns out that in case of unstable spatial unilateral
 autoregressive processes asymptotic normality holds, see Theorem \ref{main}.

With the help of the stability parameter \ $\varrho$  \  the model
can also be  written in the form
 \begin{equation}
    \label{rhomodel}
   X_{k,\ell} = \begin{cases}
               \alpha \big(X_{k-1, \ell} - \Sign(\alpha\beta)X_{k,
                 \ell-1}\big) + \varrho \Sign(\beta)X_{k, \ell-1} 
                + \varepsilon_{k, \ell} ,
                & \text{for \ $k + \ell \geq 1$,}\\
               0,  & \text{for \ $k + \ell = 0$.}
              \end{cases}
 \end{equation}
This reparametrization can be called the canonical form of \citet{ssw}
(see also \citet[17.7.6]{hamilton}). Observe that \eqref{rhomodel}
gives four different models according to the signs of \ $\alpha$ \ and
\ $\beta$. \ Hence, in order to derive estimators of the parameters \
$(\alpha,\varrho)$ \ one should have information about these signs.

For a set \ $H \subset \{(k, \ell) \in {\mathbb Z}^2 : k + \ell \geq
1\}$, \ the least 
 squares estimator \ $(\widehat\alpha_H, \widehat\varrho_H)$ \ of \
 $(\alpha, \varrho)$ 
 \ based on the observations \ $\{X_{k,\ell} : (k, \ell) \in H\}$ \ can be
 obtained by minimizing the sum of squares
 \begin{equation*}
   \sum_{(k, \ell) \in H}
    \Bigl[ X_{k, \ell} - \alpha \big(X_{k-1, \ell} -
    \Sign(\alpha\beta)X_{k, \ell-1}\big) 
           - \varrho \Sign(\beta)X_{k, \ell-1} \Bigr]^2
 \end{equation*}
 with respect to \ $\alpha$ \ and \ $\varrho$, \ and it has the form
 \begin{equation*}
   \begin{bmatrix} \widehat\alpha_H \\ \widehat\varrho_H \end{bmatrix}
   = A_H^{-1} b_H ,
 \end{equation*}
 where
 \begin{align*}
  A_H :=& \sum_{(k, \ell) \in H}
          \begin{bmatrix}
           X_{k-1, \ell} - \Sign(\alpha\beta)X_{k, \ell-1}      
             \\
            \Sign(\beta)X_{k, \ell-1}
          \end{bmatrix}
          \begin{bmatrix}
           X_{k-1, \ell} - \Sign(\alpha\beta)X_{k, \ell-1}      
             \\
            \Sign(\beta)X_{k, \ell-1}
          \end{bmatrix}^{\top}= BA^*_HB^{\top},
                   \\[2mm]
  b_H :=& \sum_{(k, \ell) \in H} X_{k, \ell}
          \begin{bmatrix}
           X_{k-1, \ell} - \Sign(\alpha\beta)X_{k, \ell-1} \\
           \Sign(\beta)X_{k, \ell-1}
          \end{bmatrix}=Bb^*_H, \qquad \text{with} \qquad
B:=\begin{bmatrix}
           1& - \Sign(\alpha\beta)\\                 
           0&  \Sign(\beta)
          \end{bmatrix}.  
 \end{align*}
Obviously, this estimator is well defined if \  $\alpha\beta\ne 0$ \
and then we have
\begin{equation*}
   \begin{bmatrix} \widehat\alpha_H \\ \widehat\varrho_H \end{bmatrix}
   = \big(B^{\top}\big)^{-1}
   \begin{bmatrix} \widehat\alpha^*_H \\ \widehat\beta^*_H \end{bmatrix}.
 \end{equation*}
Now, let us define an estimator of \
 $\beta$  \ by \ $\widehat\beta_H:=\big(\widehat\varrho_H
 -\Sign(\alpha)\widehat\alpha_H\big)\Sign(\beta)$. \ Short calculation
 shows that \ $\widehat\alpha_H =\widehat\alpha^*_H$ \ and  \
 $\widehat\beta_H =\widehat\beta^*_H$.

For \ $k,\ell \in {\mathbb Z}$ \ with \ $k+\ell \geq 1$, \ consider
the triangle 
 \begin{equation*}
   T_{k,\ell}
   := \{ (i, j) \in {\mathbb Z}^2
         : \text{$i + j \geq1$, \ $i \leq k$ \ and \ $j \leq \ell$} \} .
 \end{equation*}
For simplicity, we shall write \ $T_n := T_{n, n}$ \ for \ $n \in {\mathbb N}$.

\begin{Thm} 
\label{main} \ Let \ $\{ \varepsilon_{k, \ell} : k, \ell \in {\mathbb Z}, \, k
  + \ell \geq 1 \}$ \ be 
 independent random variables with \ ${\mathsf E}(\varepsilon_{k,\ell}) = 0$,
 \ $\Var(\varepsilon_{k, \ell}) = 1$ \ and
 \ $\sup\{ {\mathsf E}(\varepsilon_{k, \ell}^4) : k, \ell \in {\mathbb
   Z}, \, k + \ell \geq 1 \} <\infty$. \ 

If \ $|\alpha| + |\beta| < 1$ \ and \ $\alpha\beta\ne 0$ \ then
 \begin{equation}
     \label{rholimst}   
    n (\widehat\varrho_{T_n} - 1)
   \distr {\mathcal N}\left(0,
     \big(1+\kappa_{\alpha,\beta}\Sign(\alpha \beta)\big)^{-1} 
                       \sigma_{\alpha,\beta}^{-2}\right) 
 \end{equation}
and
\begin{equation}
    \label{ablimst}
   \begin{bmatrix}
    n(\widehat\alpha_{T_n} - \alpha) \\
    n(\widehat\beta_{T_n} - \beta)
   \end{bmatrix}
   \distr
   {\mathcal N}\left(\begin{bmatrix} 0 \\ 0 \end{bmatrix} ,
            \frac
            1{2\sigma_{\alpha,\beta}^2(1-\kappa_{\alpha,\beta}^2)} \begin{bmatrix}
              1 & -\kappa_{\alpha,\beta} \\ 
              -\kappa_{\alpha,\beta} & 1 \end{bmatrix} 
      \right) , 
\end{equation}
as \ $n\to\infty$, \ where
\begin{align*}
  \sigma_{\alpha,\beta}^2
  &:=\big((1+\alpha+\beta)(1+\alpha-\beta)(1-\alpha+\beta)
                               (1-\alpha-\beta)\big)^{-1/2},\\[2mm]
   \kappa_{\alpha,\beta} 
   &:=\frac{(1-\alpha^2-\beta^2)\sigma_{\alpha,\beta}^2-1}
                         {2\alpha\beta\sigma_{\alpha,\beta}^2}.
 \end{align*}

If \ $|\alpha| + |\beta| = 1$ \ and \ $\alpha\beta\ne 0$ \ then
 \begin{equation}
     \label{rholim}   
    n^{5/4} (\widehat\varrho_{T_n} - 1)
   \distr {\mathcal N}(0, \psi_\alpha) \qquad
 \end{equation}
and
\begin{equation}
    \label{ablim}
   \begin{bmatrix}
    n(\widehat\alpha_{T_n} - \alpha) \\
    n(\widehat\beta_{T_n} - \beta)
   \end{bmatrix}
   \distr
   {\mathcal N}\left(\begin{bmatrix} 0 \\ 0 \end{bmatrix} ,
            \varphi_\alpha \begin{bmatrix} 1 & -\Sign(\alpha\beta) \\
              -\Sign(\alpha\beta) & 1 \end{bmatrix} 
      \right) , 
\end{equation}
as \ $n\to\infty$, \ where
 \begin{equation*}  
   \varphi_\alpha := \frac{|\alpha|(1-|\alpha|)}{2} \qquad \text{and} \qquad 
   \psi_\alpha := \frac{15 \sqrt{\pi |\alpha|(1-|\alpha|)}}{2^{9/2}} . 
 \end{equation*}
\end{Thm}

We remark that \eqref{rholimst} is a direct consequence of
\eqref{ablimst} which is the first statement of Theorem 1.1 of
\citet{bpz2}. Further, \eqref{ablim} has already proved in
\citet{bpz2}, too, but with a far more complicated method than here.

Now, observe that to prove the statement of Theorem \ref{main} in the
unstable case it suffices to show 
\begin{equation}
     \label{arholim}
 \begin{bmatrix}
    n (\widehat\alpha_{T_n} - \alpha) \\
    n^{5/4} (\widehat\varrho_{T_n} - 1)
   \end{bmatrix}
   \distr {\mathcal N}(0, \Sigma_\alpha) \qquad
   \text{as \ $n\to\infty$,}
 \end{equation}
where
 \begin{equation*}
   \Sigma_\alpha := \begin{bmatrix}
                   \varphi_\alpha & 0 \\
                   0 & \psi_\alpha
                  \end{bmatrix}.   
 \end{equation*}
Limit \eqref{rholim} is obvious, but  \eqref{arholim} also implies \
$n(\widehat\varrho_{T_n} - 
1) \stoch 0$. \ In this way \eqref{ablim} follows from
\eqref{arholim} and
 \begin{align*}
   \begin{bmatrix}
    n (\widehat\alpha_{T_n} - \alpha) \\
    n (\widehat\beta_{T_n} - \beta)
   \end{bmatrix}
   &= \begin{bmatrix}
      n (\widehat\alpha_{T_n} - \alpha) \\
      \Sign(\beta)\big(n (\widehat\varrho_{T_n} - 1) - \Sign(\alpha)n
      (\widehat\alpha_{T_n} - \alpha) \big)
     \end{bmatrix} \\
   &= n (\widehat\alpha_{T_n} - \alpha) \begin{bmatrix} 1 \\
     -\Sign(\alpha\beta) \end{bmatrix} 
     + n (\widehat\varrho_{T_n} - 1) \begin{bmatrix} 0 \\
       \Sign(\beta) \end{bmatrix}. 
 \end{align*}

Further, we can write
 \begin{equation*}
   \begin{bmatrix} \widehat\alpha_{T_n} - \alpha \\
     \widehat\varrho_{T_n} - 1 \end{bmatrix} 
   = A_{T_n}^{-1} d_{T_n} ,
 \end{equation*}
 where
 \begin{gather*}
  d_{T_n} := \sum_{(k, \ell) \in T_n}
              \varepsilon_{k, \ell}\begin{bmatrix}
              X_{k-1, \ell} -  \Sign(\alpha\beta)X_{k, \ell-1}  \\
              \Sign(\beta) X_{k, \ell-1} 
             \end{bmatrix} ,
 \end{gather*}
and by the continuous mapping theorem \eqref{arholim}
is a consequence of the convergence
 \begin{equation}
   \label{contmap} 
  \big({\widetilde A}_{T_n}, {\widetilde d}_{T_n}\big)
  := 
  \left( \begin{bmatrix} n^{-1} & 0 \\ 0 & n^{-5/4} \end{bmatrix}
         A_{T_n} \begin{bmatrix} n^{-1} & 0 \\ 0 & n^{-5/4} \end{bmatrix} ,
         \begin{bmatrix} n^{-1} & 0 \\ 0 & n^{-5/4} \end{bmatrix} d_{T_n} \right)
  \distr 
  \big ({\widetilde A}, {\widetilde d}\, \big)
 \end{equation} 
 as \ $n\to\infty$, \ where
 \begin{equation*}
   {\widetilde A} := \begin{bmatrix}
           1/\varphi_\alpha & 0 \\
           0 & 1/\psi_\alpha
          \end{bmatrix}  \qquad \text{and} \qquad
   {\widetilde d} \distre {\mathcal N}(0, {\widetilde A}) .
 \end{equation*}
Obviously, \eqref{contmap} can be verified by proving the following
two propositions.
\begin{Pro}
  \label{An}
\ Under the conditions of Theorem \ref{main} 
\begin{equation}
   \label{Anlim}
  {\widetilde A}_{T_n} \stoch {\widetilde A} \qquad \text{as \
    $n\to\infty$}.
\end{equation}
\end{Pro}
\begin{Pro}
  \label{Dn}
\ Under the conditions of Theorem \ref{main} 
\begin{equation*}
{\widetilde d}_{T_n} \distr {\mathcal N}(0, {\widetilde A}) \qquad \text{as \
    $n\to\infty$}.
\end{equation*}
\end{Pro}

The aim of the following discussion is to show that it suffices to prove
Propositions \ref{An} and \ref{Dn} for \ $\alpha >0$ \ and \ $\beta >0$ 
\ implying \ $\varrho=\alpha+\beta$. \ In this case we have
 \begin{align}
     \label{Aplus}
  {\widetilde A}_{T_n} &= \sum_{(k, \ell) \in T_n}
              \begin{bmatrix}
               n^{-2} \big(X_{k-1, \ell} -  X_{k, \ell-1}\big)^2         
                & n^{-9/4} \big(X_{k-1, \ell} -  X_{k,
                  \ell-1}\big)  X_{k, \ell-1} \\ 
               n^{-9/4} \big(X_{k-1, \ell} - X_{k,
                 \ell-1}\big)  X_{k,\ell-1} 
                & n^{-5/2} X_{k, \ell-1}^2
              \end{bmatrix} , \\[2mm]
  {\widetilde d}_{T_n} &= \sum_{(k, \ell) \in T_n}
              \begin{bmatrix}
               n^{-1} \big(X_{k-1, \ell} - X_{k,
                 \ell-1}\big) \varepsilon_{k, \ell} \\ 
               n^{-5/4} X_{k, \ell-1} \varepsilon_{k, \ell}
              \end{bmatrix} . \label{Dplus}
 \end{align}
Model equation \eqref{model} implies that
random variable \ $X_{k,\ell}$ \ can be expressed as 
a linear combination of the variables \ $\{\varepsilon_{i,j}:(i,j)\in
T_{k,\ell}\}$, \ namely,
 \begin{equation}
     \label{MArep}
  X_{k,\ell}=\sum_{(i,j)\in T_{k,\ell}}\binom{k+\ell-i-j}{k-i}
                     \alpha^{k-i}\beta^{\ell-j}\varepsilon_{i,j}
 \end{equation}
for \ $k,\ell\in{\mathbb Z}$ \ with \ $k+\ell\geq1$. \ If \ $\alpha+\beta=1$ \
we can also write 
\begin{equation}
     \label{MArepprob}
  X_{k,\ell}=\sum_{(i,j)\in T_{k,\ell}}{\mathsf P}
  \big(S^{(\alpha)}_{k+\ell -i-j}=k-i\big)\varepsilon_{i,j}, 
 \end{equation}
where \ $S^{(\alpha)}_n$ \ is a binomial random variable with
parameters \ $(n,\alpha)$. 

Let \ $\alpha <0, \ \beta <0$ \ implying \ $\varrho
=-\alpha-\beta$ \ and  put \
$\varepsilon_{k,\ell}^{*}:=(-1)^{k+\ell}\varepsilon_{k,\ell}$ \
for  \ $k,\ell\in{\mathbb Z}$ \ with \ $k+\ell\geq1$.
\ Then \ $\{\varepsilon_{k,\ell}^{*} :k,\ell\in{\mathbb
  Z},\,k+\ell\geq1\}$ \ are independent  
random variables with \ ${\mathsf E}\,(\varepsilon_{k,\ell}^{*})=0$, \
and \ $\Var\,(\varepsilon_{k,\ell}^{*})=1$.  
\ Consider the zero start triangular spatial AR process
  \ $\{X_{k,\ell}^{*}:k,\ell\in{\mathbb Z},\,k+\ell\geq0\}$ \
  defined by
\begin{equation*}
 X_{k,\ell}^{*}
   =\begin{cases}
     -\alpha X_{k-1,\ell}^{*}-\beta
     X_{k,\ell-1}^{*}+\varepsilon_{k,\ell}^{*},
      & \text{for \ $k+\ell\geq1$,}\\
     0,  & \text{for \ $k+\ell=0$.}
    \end{cases}
\end{equation*}
In this case \eqref{rhomodel} takes the form
\begin{equation*}
   X_{k,\ell}^{*} = \begin{cases}
               -\alpha \big(X_{k-1, \ell}^{*} - X_{k,
                 \ell-1}^{*}\big) + \varrho X_{k, \ell-1}^{*} 
                + \varepsilon_{k, \ell}^{*} ,
                & \text{for \ $k + \ell \geq 1$,}\\
               0,  & \text{for \ $k + \ell = 0$.}
              \end{cases}
 \end{equation*}
Then, by representation \eqref{MArep}, 
\begin{equation*}
X_{k,\ell}^{*}
   =\sum_{(i,j)\in T_{k,\ell}}\binom{k+\ell-i-j}{k-i}
    (-\alpha)^{k-i}(-\beta)^{\ell-j}\varepsilon_{i,j}^{*}
   =(-1)^{k+\ell}X_{k,\ell}
\end{equation*}
for \ $k,\ell\in{\mathbb Z}$ \ with \ $k+\ell\geq0$. \
\ Hence,
\begin{align*}
  A_{T_n}^{*} &= \sum_{(k, \ell) \in T_n}
              \begin{bmatrix}
               n^{-2} \big(X_{k-1, \ell}^{*} - X_{k, \ell-1}^{*}\big)^2         
                & -n^{-9/4} \big(X_{k-1, \ell}^{*} - X_{k,
                  \ell-1}^{*}\big)  X_{k, \ell-1}^{*} \\ 
               -n^{-9/4} \big(X_{k-1, \ell}^{*} -  X_{k,
                 \ell-1}^{*}\big)  X_{k,\ell-1}^{*} 
                & n^{-5/2} \big(X_{k, \ell-1}^{*}\big)^2
              \end{bmatrix} \\ 
              &= 
                \begin{bmatrix}
                 -1 & 1 \\
                  1 & -1  
                 \end{bmatrix}  \widetilde A_{T_n}
                 \begin{bmatrix}
                 -1 & 1 \\
                  1 & -1  
                 \end{bmatrix}, \\[2mm]
  d_{T_n}^{*} &= \sum_{(k, \ell) \in T_n}
              \begin{bmatrix}
               -n^{-1} \big (X_{k-1, \ell}^{*} -  X_{k,
                 \ell-1}^{*}\big) \varepsilon_{k, \ell}^{*} \\ 
               n^{-5/4}  X_{k, \ell-1}^{*} \varepsilon_{k, \ell}^{*}
              \end{bmatrix} =\begin{bmatrix}
                 -1 & 1 \\
                  1 & -1  
                 \end{bmatrix}  \widetilde d_{T_n},
 \end{align*}
where \ $ \widetilde A_{T_n}$ \ and \ $ \widetilde d_{T_n}$ \ have
forms \eqref{Aplus} and \eqref{Dplus}, respectively.
Consequently, in order to prove Propositions \ref{An} and \ref{Dn} for
\ $\alpha <0$ \ and \ $\beta <0$  \ it 
suffices to prove them for \ $\alpha>0$ \ and \ $\beta>0$.

Next, let \ $\alpha <0, \ \beta >0$ \ implying \ $\varrho
=-\alpha+\beta$ \ and  put \
$\varepsilon_{k,\ell}^{+}:=(-1)^k\varepsilon_{k,\ell}$ \
for  \ $k,\ell\in{\mathbb Z}$ \ with \ $k+\ell\geq1$.
\ Then \ $\{\varepsilon_{k,\ell}^{+} :k,\ell\in{\mathbb
  Z},\,k+\ell\geq1\}$ \ are again independent  
random variables with \ ${\mathsf E}\,(\varepsilon_{k,\ell}^{+})=0$, \
and \ $\Var\,(\varepsilon_{k,\ell}^{+})=1$.  
\ Consider the zero start triangular spatial AR process
  \ $\{X_{k,\ell}^{+}:k,\ell\in{\mathbb Z},\,k+\ell\geq0\}$ \
  defined by
\begin{equation*}
 X_{k,\ell}^{+}
   =\begin{cases}
     -\alpha X_{k-1,\ell}^{+}+\beta
     X_{k,\ell-1}^{+}+\varepsilon_{k,\ell}^{+},
      & \text{for \ $k+\ell\geq1$,}\\
     0,  & \text{for \ $k+\ell=0$.}
    \end{cases}
\end{equation*}
Now \  $X_{k,\ell}^{+}=(-1)^k X_{k,\ell}$ \ and \eqref{rhomodel} takes the form
\begin{equation*}
   X_{k,\ell}^{+} = \begin{cases}
               -\alpha \big(X_{k-1, \ell}^{+} - X_{k,
                 \ell-1}^{+}\big) + \varrho X_{k, \ell-1}^{+} 
                + \varepsilon_{k, \ell}^{+} ,
                & \text{for \ $k + \ell \geq 1$,}\\
               0,  & \text{for \ $k + \ell = 0$}
              \end{cases}
 \end{equation*}
and
\begin{align*}
  A_{T_n}^{+} &= \sum_{(k, \ell) \in T_n}
              \begin{bmatrix}
               n^{-2} \big(X_{k-1, \ell}^{+} + X_{k, \ell-1}^{+}\big)^2         
                & -n^{-9/4} \big(X_{k-1, \ell}^{+} + X_{k,
                  \ell-1}^{+}\big)  X_{k, \ell-1}^{+} \\ 
               -n^{-9/4} \big(X_{k-1, \ell}^{+} +  X_{k,
                 \ell-1}^{+}\big)  X_{k,\ell-1}^{+} 
                & n^{-5/2} \big(X_{k, \ell-1}^{+}\big)^2
              \end{bmatrix} \\ 
              &= 
                \begin{bmatrix}
                 -1 & 1 \\
                  1 & -1  
                 \end{bmatrix}  \widetilde A_{T_n}
                 \begin{bmatrix}
                 -1 & 1 \\
                  1 & -1  
                 \end{bmatrix}, \\[2mm]
  d_{T_n}^{+} &= \sum_{(k, \ell) \in T_n}
              \begin{bmatrix}
               -n^{-1} \big (X_{k-1, \ell}^{+} +  X_{k,
                 \ell-1}^{+}\big) \varepsilon_{k, \ell}^{+} \\ 
               n^{-5/4}  X_{k, \ell-1}^{+} \varepsilon_{k, \ell}^{+}
              \end{bmatrix} =\begin{bmatrix}
                 -1 & 1 \\
                  1 & -1  
                 \end{bmatrix}  \widetilde d_{T_n},
 \end{align*}
where \ $ \widetilde A_{T_n}$ \ and \ $ \widetilde d_{T_n}$ \ have
forms \eqref{Aplus} and \eqref{Dplus}, respectively.

In the same way one can handle the case \ $\alpha >0, \ \ \beta <0$ \
implying \ $\varrho=\alpha-\beta$ \ by considering \
$\{X_{k,\ell}^{\circ}:k,\ell\in{\mathbb Z},\,k+\ell\geq0\}$ \ 
  defined by
\begin{equation*}
 X_{k,\ell}^{\circ}
   =\begin{cases}
     \alpha X_{k-1,\ell}^{\circ}-\beta
     X_{k,\ell-1}^{\circ}+\varepsilon_{k,\ell}^{\circ},
      & \text{for \ $k+\ell\geq1$,}\\
     0,  & \text{for \ $k+\ell=0$,}
    \end{cases}
\end{equation*}
with  $\varepsilon_{k,\ell}^{\circ}=(-1)^{\ell}
\varepsilon_{k,\ell}$. 

\section{Results on the covariance structure}
  \label{sec:sec2}
In order to prove Propositions \ref{An} and \ref{Dn} one has to know
the asymptotic behaviour of the covariances of the process \
$X_{k,\ell}$. \ By representation \eqref{MArep} we obtain that for all
\ $k_1,\ell_1,k_2,\ell_2\in{\mathbb Z}$ \ with \ $k_1+\ell_1\geq0$ \ and
\ $k_2+\ell_2\geq0$, \ and for all \ $\alpha,\beta\in {\mathbb R}$, 
\begin{equation}
  \label{eq:eq2.1}
  \Cov\big(X_{k_1,\ell_1},X_{k_2,\ell_2}\big)
 =\!\!\!\!\!\! \sum_{(i,j)\in T_{k_1\land k_2,\ell_1\land \ell_2}}\!\!
        \binom{k_1+\ell_1\!-\!i\!-\!j}{k_1\!-\!i}
        \binom{k_2\!+\!\ell_2\!-\!i\!-\!j}{k_2\!-\!i} 
        \alpha^{k_1+k_2-2i}\beta^{\ell_1+\ell_2-2j},
\end{equation}
where \ $k\land \ell :=\min\{k,\ell\}$ \ and an empty sum is defined
to be equal to $0$. \ Observe,  if \
$0<\alpha <1$ \ and \ $\beta=1-\alpha$ \ then by representation
\eqref{MArepprob} covariance \eqref{eq:eq2.1} can be
expressed in the form
\begin{equation*}
\Cov\big(X_{k_1,\ell_1},X_{k_2,\ell_2}\big)=\sum_{m=1}^{k_1\land
  k_2+\ell_1\land \ell_2} {\mathsf
  P}\big(S_{k_1+\ell_1-m,k_2+\ell_2-m}^{(\alpha,1-\alpha)}=k_1+\ell_2-m\big), 
\end{equation*}
where for \  $\nu,\mu \in (0,1)$ \  real numbers  \
$S_{k,\ell}^{(\mu,\nu)}:=\xi_k^{(\mu)}+\eta_{\ell }^{(\nu)}$, \ 
and \ $\xi_k^{(\mu)}$ \ and  \ 
$\eta_{\ell }^{(\nu)}$ \ are independent binomial random variables
with parameters \ $(k,\mu)$ \ and \ $(\ell,\nu)$, \
respectively. Now, Lemmas 2.4 and 2.6 of \citet{bp1} directly imply
that there exists a constant \
$D_{\mu,\nu}>0$ \ such that for all \ 
$k,\ell \geq 0, \ k+\ell\geq 1$, \ $0\leq i\leq k+\ell$ \ and \ $0\leq
j\leq k+\ell -1$ we have 
\begin{equation}
  \label{eq:eq2.2}
{\mathsf P}\big(S_{k,\ell}^{(\mu,\nu)}=i\big) \leq \frac
{D_{\mu,\nu}}{\sqrt{k+\ell}}
\qquad \text{and} \qquad
\Big|{\mathsf P}\big(S_{k,\ell}^{(\mu,\nu)}=j+1\big)-{\mathsf
  P}\big(S_{k,\ell}^{(\mu,\nu)}=j\big)\Big| \leq \frac {D_{\mu,\nu}}{k+\ell}. 
\end{equation}

Hence, one can 
determine the magnitudes of the covariances and prove the following
lemma.  
\begin{Lem}
   \label{covbound}
\citep[Lemma 2.1]{bpz2} \ If \ $|\alpha|+|\beta|=1$ \ and \
$0<|\alpha|<1$ \ then 
\begin{equation*}
 \big|\Cov\big(X_{k_1,\ell_1},X_{k_2,\ell_2}\big)\big|
    \leq C_{\alpha}\sqrt{k_1+\ell_1+k_2+\ell_2}
\end{equation*}
with some constant \ $C_{\alpha}>0$.
\end{Lem}

Now, for \ $n\in{\mathbb N}$, \ let us introduce piecewise
constant random fields 
 \begin{equation*}
  Z_{1,0}^{(n)}(s,t):=n^{-1/4}X_{[ns]+1,[nt]}
  \quad \text{and} \quad Z_{0,1}^{(n)}(s,t):=n^{-1/4}X_{[ns],[nt]+1}, \qquad s,t\in
  {\mathbb R}, \ s+t\geq0.  
 \end{equation*} 
Concerning the asymptotic behaviour of their covariances
one can verify the following result.
\begin{Pro}
  \label{covlim}
\citep[Proposition 2.2]{bpz2} \ Let \ $s_1,t_1,s_2,t_2\in{\mathbb R}$
\ with \ $s_1+t_1>0$, \ $s_2+t_2>0$.  
\ If \ $0<\alpha<1$ \ and \ $\beta=1-\alpha$ \ then
\begin{equation*}
    \begin{bmatrix}
     \Cov(Z_{1,0}^{(n)}(s_1,t_1),Z_{1,0}^{(n)}(s_2,t_2))
     & \Cov(Z_{1,0}^{(n)}(s_1,t_1),Z_{0,1}^{(n)}(s_2,t_2)) \\
     \Cov(Z_{1,0}^{(n)}(s_2,t_2),Z_{0,1}^{(n)}(s_1,t_1))
    & \Cov(Z_{0,1}^{(n)}(s_1,t_1),Z_{0,1}^{(n)}(s_2,t_2))
   \end{bmatrix}
   \!\to\!z_\alpha(s_1,t_1,s_2,t_2) \begin{bmatrix} 1 & 0 \\ 0 & 1 \end{bmatrix}
\end{equation*}
as \ $n\to\infty$, \ where 
\begin{equation*}
 z_\alpha(s_1,t_1,s_2,t_2)
    =\begin{cases}
       \frac{\sqrt{s_1+s_2+t_1+t_2}-
                      \sqrt{|s_1-s_2|+|t_1-t_2|}}
                      {\sqrt{2\pi\alpha(1-\alpha)}}
     & \text{if \ $(1-\alpha )(s_1-s_2) \!=\!\alpha (t_1-t_2)$,}\\
     0, & \text{otherwise.}
     \end{cases}
\end{equation*}
Moreover, if \ $(1-\alpha )(s_1-s_2)\ne \alpha (t_1-t_2)$ \ then the 
convergence to \ $0$ \ has an exponential rate.
\end{Pro}

Further, one can also estimate the difference of two neighbouring
covariances.
\begin{Pro}
  \label{covdiff}
 \citep[Proposition 2.5]{bpz2}
\ If  \ $0<\alpha<1$ \ and \ $\beta=1-\alpha$ \ then there exists a constant
$K_{\alpha}>0$ such that
\begin{equation*}
\big |\Cov (Z_{i,j}^{(n)}(s_1,t_1),Z_{j,i}^{(n)}(s_2,t_2))-
\Cov (Z_{i,j}^{(n)}(s_1,t_1),Z_{i,j}^{(n)}(s_2,t_2))\big |\leq
K_{\alpha }n^{-1/2}
\end{equation*}
for all \ $n\in{\mathbb N}, \ s_1,t_1,s_2,t_2\in {\mathbb R}$, \ with 
$s_1+t_1>0,\ s_2+t_2>0$ \ and $(i,j)\in \big\{ (0,1),(1,0) \big\}$. 
\end{Pro}

Finally, in order to estimate covariances we make use of he following 
 lemma which is a generalization of \citet[Lemma 11]{bpz1}.

\begin{Lem}
   \label{mmmm}
\citep[Lemma 2.8]{bpz2}
Let \ $\xi_1,\ldots,\xi_N$ \ be independent random variables with \ 
${\mathsf E}(\xi_i)=0$,
 \ ${\mathsf E}(\xi_i^2)=1$ \ for all \ $i=1,\ldots,N$, \ and
 \ $M_4:=\max_{1\leq i\leq N}{\mathsf E}(\xi_i^4)<\infty$.
\ Let \ $a_1,\ldots,a_{n_1},b_1,\ldots,b_{n_2}$, $c_1,\ldots,c_{n_3},
d_1,\ldots,d_{n_4}\in{\mathbb R}, \ n_1,n_2,n_3,n_4\leq N$ \ and
\begin{equation*}
X:=\sum_{i=1}^{n_1}a_i\xi_i, \quad Y:=\sum_{j=1}^{n_2}b_j\xi_j, 
 \quad Z:=\sum_{i=1}^{n_3}c_i\xi_i, \quad W:=\sum_{j=1}^{n_4}d_j\xi_j.
\end{equation*}
Then
\begin{equation*}
\Cov(XY,ZW)\!=\!\!\!\!\!\!\!\!\!\!
\sum_{i=1}^{n_1\land n_2\land n_3\land n_4}\!\!\!\!\!\!\! \big({\mathsf E} 
(\xi _i^4)\!-3\big)\,  a_ib_ic_id_i+\Cov (X,Z)\Cov(Y,W)+\Cov (X,W)\Cov(Y,Z).
\end{equation*}
Moreover, if  \ $a_i,b_i,c_i,d_i\geq 0$ \ then
\begin{equation*}
0\leq\Cov (XY,ZW)\leq M_4\Cov (X,Z)\Cov(Y,W)+M_4\Cov (X,W)\Cov(Y,Z),
\end{equation*}
and
\begin{equation*}
0\leq {\mathsf E} (XYZW)\leq M_4\big({\mathsf E} (XZ)\,{\mathsf E}
(YW)+{\mathsf E} (XW)\,{\mathsf E} (YZ) +{\mathsf E} 
(XY)\,{\mathsf E} (ZW) \big).
\end{equation*}
\end{Lem}

\section{Proof of Proposition \ref{An}}
  \label{sec:sec3}
Let \ $\alpha, \beta\in (0,1)$ \ with \ $\alpha+\beta=1$ \ and
\begin{equation*}
S_{n,1}:=\!\!\!\sum_{(k, \ell) \in T_n}\!\! \big(X_{k-1, \ell} - X_{k,
  \ell-1}\big)^2, \quad 
S_{n,2}:=\!\!\!\sum_{(k, \ell) \in T_n} \!\!\big(X_{k-1, \ell} - X_{k,
  \ell-1}\big)X_{k, \ell-1}, \quad
S_{n,3}:=\!\!\!\sum_{(k, \ell) \in T_n} \!\!X_{k, \ell-1}^2.
\end{equation*}
Thus,
\begin{equation*}
\widetilde A_{T_n}=\begin{bmatrix}
                  n^{-2}S_{n,1} &  n^{-9/4}S_{n,2} \\
                   n^{-9/4}S_{n,2} &  n^{-5/2}S_{n,3}
                  \end{bmatrix}
\end{equation*}
and \eqref{Anlim} follows from
\begin{equation}
   \label{eq:eq3.1}
 n^{-2}S_{n,1} \qmean \frac 1{\varphi_{\alpha}}=\frac 2{\alpha(1-\alpha)}, \quad
 n^{-9/4}S_{n,2} \qmean 0, \quad
 n^{-2}S_{n,3} \qmean \frac 1{\psi_{\alpha}}=\frac
 {2^{9/2}}{15\sqrt{\pi\alpha(1-\alpha)}}. 
\end{equation}

The last two statements of \eqref{eq:eq3.1} have already been proved, see
\citet[Proposition 1.2 and Section 6, pp. 40-41]{bpz2}. In order to
verify the remaining 
statement one has to show
\begin{equation}
  \label{eq:eq3.2}
\lim_{n\to\infty}  n^{-2}{\mathsf E}\big(S_{n,1}\big)=\frac 1{\varphi_{\alpha}}
\qquad \text{and} \qquad
\lim_{n\to\infty}  n^{-4} \Var \big(S_{n,1}\big)=0. 
\end{equation}
It is easy to see that 
\begin{equation*}
n^{-2}{\mathsf E}\big(S_{n,1}\big) = \iint\limits_T \sqrt{n}\Big(
\Var\big(Z_{0,1}^{(n)}(s,t)\big)+\Var\big(Z_{1,0}^{(n)}(s,t)\big)-
2\Cov \big(Z_{0,1}^{(n)}(s,t),Z_{1,0}^{(n)}(s,t)\big)\Big)\,{\mathrm
  d}s\,{\mathrm d}t, 
\end{equation*}
where \ $T:=\{(s,t)\in{\mathbb R}^2:s+t\geq0,\,s\leq1,\,t\leq1\}$, \
and using \eqref{eq:eq2.1} one can prove 
\begin{equation*}
\lim_{n\to\infty} \sqrt{n}\Big(\Var\big(Z_{i,j}^{(n)}(s,t)\big)- \Cov
\big(Z_{i,j}^{(n)}(s,t),Z_{j,i}^{(n)}(s,t)\big)\Big)=\frac 1{2\alpha(1-\alpha)},
\end{equation*}
where \ $(i,j)\in \big\{(0,1),(1,0)\big\}$. \ The details can
be found in \citet[Section 5, pp. 36-37]{bpz2}.
Hence, Proposition \ref{covdiff} and the dominated convergence theorem
imply the first statement of \eqref{eq:eq3.2}.

Now, by Lemma  \ref{mmmm}
\begin{align*}
\Var \big(S_{n,1}\big)&=\!\!\sum_{(k_1,\ell_1)\in
  T_n}\sum_{(k_2,\ell_2)\in T_n} \Cov \big((X_{k_1-1, \ell_1} -
X_{k_1,\ell_1-1})^2, (X_{k_2-1, \ell_2} -X_{k_2,\ell_2-1})^2\big)\\
&\leq \sum_{(k_1,\ell_1)\in
  T_n}\sum_{(k_2,\ell_2)\in T_n}\Big( 2 M_4
L^{(1)}_{k_1,\ell_1,k_2,\ell_2}+(M_4-3)^+L^{(2)}_{k_1,\ell_1,k_2,\ell_2}\Big)+{\mathcal
  O}(n^3),  
\end{align*}
where
\begin{align}
L^{(1)}_{k_1,\ell_1,k_2,\ell_2}:=&\,\Cov \big(X_{k_1-1, \ell_1} -
X_{k_1,\ell_1-1}, X_{k_2-1, \ell_2} -X_{k_2,\ell_2-1}\big)^2,
\nonumber \\
L^{(2)}_{k_1,\ell_1,k_2,\ell _2}:=&\!\!\!\!\!\!\!\!
\sum_{(i,j)\in T_{k_1\land k_2-1,\ell_1\land
    \ell_2-1}}\!\!\!\!\!\!\!\!\!\Big (
{\mathsf P}\big (S^{(\alpha )}_{k_1+\ell_1-1-i-j}=k_1-i\big )
-{\mathsf P}\big (S^{(\alpha )}_{k_1+\ell_1-1-i-j}=k_1-1-i
\big )\Big )^2  \label{eq:eq3.3} \\
&\phantom{\qquad \qquad \quad\quad  } \times\Big (
{\mathsf P}\big (S^{(\alpha )}_{k_2+\ell_2-1-i-j}=k_2-i\big )
-{\mathsf P}\big (S^{(\alpha )}_{k_2+\ell_2-1-i-j}=
k_2-1-i\big )\Big )^2  \nonumber \\
\leq&\!\!\!\!\!\!\!\!
\sum_{(i,j)\in T_{k_1\land k_2-1,\ell_1\land
    \ell_2-1}}\!\!\!\!\!\!\!\!\!\Big (
{\mathsf P}\big (S^{(\alpha )}_{k_1+\ell_1-1-i-j}=k_1-i\big )^2
+{\mathsf P}\big (S^{(\alpha )}_{k_1+\ell_1-1-i-j}=k_1-1-i
\big )^2\Big ) \nonumber\\
&\phantom{\qquad \qquad \quad\quad  } \times\Big (
{\mathsf P}\big (S^{(1-\alpha )}_{k_2+\ell_2-1-i-j}=\ell_2-1-j\big )^2
+{\mathsf P}\big (S^{(1-\alpha )}_{k_2+\ell_2-1-i-j}=
\ell_2-j\big )^2\Big ). \nonumber
\end{align}
Obviously, 
\begin{align*}
n^{-4}&\!\!\sum_{(k_1,\ell_1)\in  T_n}\sum_{(k_2,\ell_2)\in T_n}
L^{(1)}_{k_1,\ell_1,k_2,\ell_2} \\ &= \!\iint\limits_T\!\iint\limits_T\!
\Big(\sqrt{n} \Cov\big(Z_{0,1}^{(n)}(s_1,t_1)-Z_{1,0}^{(n)}(s_1,t_1),
Z_{0,1}^{(n)}(s_2,t_2)-Z_{1,0}^{(n)}(s_2,t_2)\big)\Big)^2{\mathrm d}
s_1 \, {\mathrm d} t_1 \,{\mathrm d} s_2\, {\mathrm d} t_2, 
\end{align*}
where due to Propositions \ref{covlim}, \ref{covdiff} and dominated
convergence theorem the right hand side converges to \ $0$ \ as \
$n\to\infty$. 

Further, the second inequality of \eqref{eq:eq2.2} implies 
\begin{align*}
L^{(2)}_{k_1,\ell_1,k_2,\ell_2} \leq & \sum_{(i,j)\in T_{k_1\land k_2-1,\ell_1\land
    \ell_2-1}} \frac
{D_{\alpha,\alpha}^4}{(k_1+\ell_1-1-i-j)^2(k_2+\ell_2-1-i-j)^2} \\
\leq & \sum_{m=1}^{k_1\land k_2+\ell_1\land
    \ell_2-2} \frac {D_{\alpha,\alpha}^4(k_1\land k_2+\ell_1\land
    \ell_2-1-m)}{(k_1+\ell_1-1-m)^2(k_2+\ell_2-1-m)^2} \\
\leq & \sum_{m=1}^{k_1\land k_2+\ell_1\land
    \ell_2-2} \!\!\!\!\frac {D_{\alpha,\alpha}^4}{(k_1\land k_2+\ell_1\land
    \ell_2-1-m)^3} <\!\!\!\sum_{m=1}^{k_1\land k_2+\ell_1\land
    \ell_2-2} \!\frac {D_{\alpha,\alpha}^4}{m^2}<
   \frac {\pi ^2D_{\alpha,\alpha}^4}6<\infty, 
\end{align*}
so
\begin{equation*}
n^{-4}\sum_{(k_1,\ell_1)\in  T_n}\sum_{(k_2,\ell_2)\in T_n}
L^{(2)}_{k_1,\ell_1,k_2,\ell_2} = \!\iint\limits_T\!\iint\limits_T\!
L^{(2)}_{[ns_1],[nt_1],[ns_2],[nt_2]}{\mathrm d}
s_1 \, {\mathrm d} t_1 \,{\mathrm d} s_2\, {\mathrm d} t_2 \leq \frac
{2\pi ^2D_{\alpha,\alpha}^4}3.
\end{equation*}

Finally, e.g.
\begin{align*}
\sum_{(i,j)\in T_{[ns_1]\land [ns_2]-1,[nt_1]\land
    [nt_2]-1}}\!\!\!\!\!\!\!\!\!
&{\mathsf P}\big (S^{(\alpha )}_{[ns_1]+[nt_1]-1-i-j}=[ns_1]-i\big )^2
{\mathsf P}\big (S^{(1-\alpha )}_{[ns_2]+[nt_2]-1-i-j}=[nt_2]-1-j
\big )^2 \\
&\leq \sqrt{n} \Cov\big(Z_{1,0}^{(n)}(s_1,t_1)-Z_{1,0}^{(n)}(s_2,t_2)\big),
\end{align*}
which by Proposition \ref{covlim} converges to \ $0$ \ as \
$n\to\infty$ \ if \ $(1-\alpha)(s_1-s_2)\ne \alpha (t_1-t_2)$. \
Similar results can be derived for the remaining three terms of the
right hand side of \eqref{eq:eq3.3}, so by the dominated convergence
theorem 
\begin{equation}
   \label{eq:eq3.4}
\lim_{n\to\infty} n^{-4}\sum_{(k_1,\ell_1)\in  T_n}\sum_{(k_2,\ell_2)\in T_n}
L^{(2)}_{k_1,\ell_1,k_2,\ell_2} =0,
\end{equation}
which completes the proof. \proofend

\section{Proof of Proposition \ref{Dn}}
  \label{sec:sec4}
Again, let \ $\alpha, \beta \in (0,1)$ \ with \ $\alpha+\beta=1$ \ and
denote by \ $d_n^{(i)}, \ i=1,2,$ \ the components of  \ $d_{T_n}$. \
First we show that \ $(d_{T_n})_{n\geq1}$ \ is a square integrable  
two dimensional martingale  with respect to filtration \
$({\mathcal F}_n)_{n\geq1}$, \ where \ ${\mathcal F}_n$ 
 \ denotes the \ $\sigma$--algebra generated by random variables
 \ $\{\varepsilon_{k,\ell}:(k,\ell)\in T_n\}$. 

In order to do this we give a useful decomposition of  \ $d_{T_n}-d_{T_{n-1}}$, \ 
where \ $d_{T_0}:=(0,0)^{\top}$. \ By representation \eqref{MArep},
\begin{align*}
  d_n^{(1)}-d_{n-1}^{(1)}
  =&\sum_{(k,\ell)\in T_n\setminus T_{n-1}}\varepsilon_{k,\ell}
               \Bigg( \sum_{(i,j)\in T_{k-1,\ell}}
                {\mathsf P}\big (S^{(\alpha )}_{k+\ell-1-i-j}=k-1-i\big )             
                 \varepsilon_{i,j}\\
  &\phantom{==========}
                - \sum_{(i,j)\in T_{k,\ell-1}}
                {\mathsf P}\big (S^{(\alpha )}_{k+\ell-1-i-j}=k-i\big )             
                 \varepsilon_{i,j}\Bigg), \\
 d_n^{(2)}-d_{n-1}^{(2)}
  =&\sum_{(k,\ell)\in T_n\setminus T_{n-1}}\varepsilon_{k,\ell}
                \sum_{(i,j)\in T_{k,\ell -1}}                
                 {\mathsf P}\big (S^{(\alpha )}_{k+\ell-1-i-j}=k-i\big
                 ) \varepsilon_{i,j}.   
 \end{align*}
Collecting first the terms containing only \ $\varepsilon_{i,j}$ \ with
 \ $(i,j)\in T_n\setminus T_{n-1}$, \ and then the rest, we obtain 
 decomposition 
\begin{equation}
       \label{eq:eq4.1}
  d_{T_n}-d_{T_{n-1}}
  =d_{n,1}+\sum_{(k,\ell)\in T_n\setminus T_{n-1}}
            \varepsilon_{k,\ell} \,d_{n,2,k,\ell},  
 \end{equation}
 where \
 $d_{n,1}=\big(\delta_{n,1}^{(1)}-\delta_{n,1}^{(2)},\delta_{n,1}^{(2)}\big)^{\top}$
 \ and \ 
$d_{n,2,k,\ell}=\big (\delta_{n,2,k-1,\ell}-\delta_{n,2,k,\ell-1 },
\delta_{n,2,k,\ell-1 }\big )^{\top}$ \ with  
\begin{align*}
  \delta_{n,1}^{(1)}
  &:=\!\!\!\!\sum_{(k,\ell)\in T_n\setminus T_{n-1}}\!\!\!\!\varepsilon_{k,\ell}\!\!\!\!
                \sum_{(i,j)\in T_{k-1,\ell}\setminus T_{n-1}}\!\!\!\!
                 {\mathsf P}\big (S^{(\alpha )}_{k+\ell-1-i-j}=k\!-\!1\!-\!i\big )  
                 \varepsilon_{i,j} 
 =\!\! \sum_{k=-n+2}^m\sum_{i=-n+1}^{k-1}\!\!\alpha ^{k-1-i}\varepsilon
 _{k,n}\varepsilon _{i,n}, \\
 \delta_{n,1}^{(2)}
  &:=\!\!\!\!\sum_{(k,\ell)\in T_n\setminus T_{n-1}}\!\!\!\!\varepsilon_{k,\ell}\!\!\!\!
                \sum_{(i,j)\in T_{k,\ell -1}\setminus T_{n-1}}\!\!\!\!
                {\mathsf P}\big (S^{(\alpha )}_{k+\ell-1-i-j}=k\!-\!i\big )  
                \varepsilon_{i,j}  
 =\!\!\sum_{\ell =-n+2}^n\sum_{j=-n+1}^{\ell -1}\!\!\beta ^{\ell-1-j}
    \varepsilon _{n,\ell}\varepsilon _{n,j}, \\
  \delta_{n,2,k,\ell}
  &:=\sum_{(i,j)\in T_{k,\ell}\cap\,T_{n-1}} 
                  {\mathsf P}\big (S^{(\alpha )}_{k+\ell-i-j}=k-i\big )
                 \varepsilon_{i,j}. 
 \end{align*}
The components of  \ $d_{n,1}$ \ are quadratic forms of the variables
 \ $\{\varepsilon_{i,j}:(i,j)\in T_n\setminus T_{n-1}\}$, \ hence \ 
$d_{n,1}$ \ is independent of \ ${\mathcal F}_{n-1}$.
\ Besides this the terms \ $\delta_{n,2,k,\ell}$ \ are linear 
combinations of the variables
 \ $\{\varepsilon_{i,j}:(i,j)\in T_{n-1}\}$, \ thus 
\ vectors \ $d_{n,2,k,\ell}$  are measurable with respect 
to \ ${\mathcal F}_{n-1}$. \ Consequently,
\begin{equation*}
{\mathsf E}(d_{T_n}-d_{T_{n-1}} \mid {\mathcal F}_{n-1})
  ={\mathsf E} (d_{n,1})
   +\sum_{(k,\ell)\in T_n\setminus T_{n-1}}
     d_{n,2,k,\ell}\,{\mathsf E}(\varepsilon_{k,\ell}  \mid {\mathcal F}_{n-1})
  =0.
\end{equation*}
Hence \ $(d_{T_n})_{n\geq1}$ \ is a square integrable martingale with
 respect to the filtration \ $({\mathcal F}_n)_{n\geq1}$ \ and
 obviously the same is valid for \  $(\widetilde d_{T_n})_{n\geq1}$. 

By the Martingale Central Limit Theorem \citep{jc}, in order to prove the
 statement of Proposition \ref{Dn}, it suffices to show that the conditional
 variances of the martingale differences converge in probability and to verify
 the conditional Lindeberg condition.
To be precise, the statement is a consequence of the following two
 propositions, where \ $\bone_H$ \ denotes the indicator function of a set
 \ $H$. 

\begin{Pro}
     \label{CCDn}
\begin{equation*}
\sum_{m=1}^n{\mathsf E}\Big (\big(\widetilde d_{T_m}-\widetilde d_{T_{m-1}}\big)
\big(\widetilde d_{T_m}-\widetilde d_{T_{m-1}}\big)^{\top} \,\Big |\,
{\mathcal F}_{m-1}\Big) \stoch \widetilde A \qquad\qquad \text{as \
  $n\to\infty$.} 
\end{equation*}

\end{Pro}

\smallskip
\begin{Pro}
    \label{LINDDn} \
For all \ $\delta>0$,
\begin{equation*}
\sum_{m=1}^n{\mathsf E}\Big (\big\Vert\widetilde d_{T_m} -\widetilde
d_{T_{m-1}}\big\Vert ^2  \bone_{\left\{\Vert \widetilde
    d_{T_m}-\widetilde d_{T_{m-1}}\Vert\geq \delta \right\}} \,\Big |\, 
{\mathcal F}_{m-1}\Big) \stoch 0 \qquad \qquad \text{as \ $n\to\infty$.}
\end{equation*}
\end{Pro}

\bigskip
\noindent
{\bf Proof of Proposition \ref{CCDn}.} \ 
Considering separately the entries of \ $\big(\widetilde
d_{T_m}-\widetilde d_{T_{m-1}}\big) 
\big(\widetilde d_{T_m}-\widetilde d_{T_{m-1}}\big)^{\top}$ \ one can see
that the statement of the proposition is a consequence of
\begin{align}
n^{-2}&\sum_{m=1}^n{\mathsf E}\Big (\big
(d_m^{(1)}-d_{m-1}^{(1)}\big)^2 \,\Big |\, {\mathcal F}_{m-1}\Big)
\qmean  \frac 1{\varphi_{\alpha}}, \label{eq:eq4.2} \\
n^{-5/2}&\sum_{m=1}^n{\mathsf E}\Big (\big
(d_m^{(2)}-d_{m-1}^{(2)}\big)^2 \,\Big |\, {\mathcal F}_{m-1}\Big)
\qmean  \frac 1{\psi_{\alpha}}, \label{eq:eq4.3}\\
n^{-9/4}&\sum_{m=1}^n{\mathsf E}\Big (\big
(d_m^{(1)}-d_{m-1}^{(1)}\big)\big(d_m^{(2)}-d_{m-1}^{(2)}\big) \,\Big |\,
{\mathcal F}_{m-1}\Big) \qmean  0 \label{eq:eq4.4}
\end{align}
as \ $n\to\infty$. \ Limits \eqref{eq:eq4.2} and \eqref{eq:eq4.3} have
already been proved, see \citet[Section 6, pp. 40-41 and Proposition
4.1]{bpz2}. A more detailed proof can be found in \citet[Propositions
6.1 and 4.1]{bpz3}. 

Now, let \ $U_m:= {\mathsf E}\big (
(d_m^{(1)}-d_{m-1}^{(1)})(d_m^{(2)}-d_{m-1}^{(2)}) \,\big |\,
{\mathcal F}_{m-1}\big)$  \ and we have
\begin{equation*}
d_m^{(1)}-d_{m-1}^{(1)}=\!\!\!\sum_{(k, \ell) \in T_m\setminus T_{m-1}}
(X_{k-1, \ell}-X_{k,\ell-1}) \varepsilon_{k, \ell} , \qquad
d_m^{(2)}-d_{m-1}^{(2)}=\!\!\!\sum_{(k, \ell) \in T_m\setminus T_{m-1}}
X_{k,\ell-1} \varepsilon_{k, \ell}.  
\end{equation*} 
Representation \eqref{MArep} and independence of the
error terms \ $\varepsilon_{i,j} $ \ imply 
\begin{align*}
{\mathsf E}\Big(\big
(d_m^{(1)}-d_{m-1}^{(1)}\big)&\big(d_m^{(2)}-d_{m-1}^{(2)}\big)\Big)=
\sum_{(k, \ell) \in T_m\setminus T_{m-1}}  {\mathsf E}\big((X_{k-1,
  \ell}-X_{k,\ell-1}) X_{k,\ell-1}\big) \, {\mathsf E} \big(\varepsilon_{k,
  \ell}^2\big) \\
&=
\sum_{(k, \ell) \in T_m\setminus T_{m-1}}  {\mathsf E}\big((X_{k-1,
  \ell}-X_{k,\ell-1}) X_{k,\ell-1}\big) ={\mathsf E}(S_{m,2})-{\mathsf E}(S_{m-1,2}),
\end{align*}
so using the second statement of \eqref{eq:eq3.1} we obtain
\begin{equation*}
n^{-9/4}\sum_{m=1}^n {\mathsf E}(U_m) =n^{-9/4} {\mathsf E} (S_{n,2})\to 0
\qquad\qquad \text{as \ $n\to\infty$}.
\end{equation*}

Further, decomposition \eqref{eq:eq4.1}, independence of \
$\delta_{m,1}^{(1)}, \ \delta_{m,1}^{(2)}$ \ and \ $\big\{\varepsilon _{k,\ell}, \
    (k,\ell)\in   T_m\setminus T_{m-1}\big\}$  \ from \ ${\mathcal
      F}_{m-1}$, \ and measurability of \ $\delta _{m,2,k,\ell}$ \ with
    respect to \ ${\mathcal F}_{m-1}$ \ imply
\begin{equation*}
U_m={\mathsf E}\Big(\big(\delta_{m,1}^{(1)}- \delta_{m,1}^{(2)}\big)
\delta_{m,1}^{(2)}\Big) +
\sum_{(k, \ell) \in T_m\setminus T_{m-1}}\big(\delta _{m,2,k-1,\ell}
-\delta _{m,2,k,\ell-1}\big) \delta _{m,2,k,\ell-1}. 
\end{equation*}
In this way, to complete the proof of \eqref{eq:eq4.4} one has to show
\begin{equation}
   \label{eq:eq4.5}
n^{-9/2} \Var \bigg(\sum_{m=1}^n U_m\bigg)= n^{-9/2} \Var
\bigg(\sum_{m=1}^n \sum_{(k, \ell) \in 
  T_m\setminus T_{m-1}}\!\!\!\!\!(\delta _{m,2,k-1,\ell} 
-\delta _{m,2,k,\ell-1}) \delta _{m,2,k,\ell-1}\bigg)\!\to\! 0
\end{equation}
as \ $n\to\infty$. 

Now, consider
\begin{align}  
\Var& \bigg(\sum_{m=1}^n \sum_{(k, \ell) \in
  T_m\setminus T_{m-1}}(\delta _{m,2,k-1,\ell} 
-\delta _{m,2,k,\ell-1}) \delta _{m,2,k,\ell-1}\bigg) \nonumber\\
&= \sum_{m_1=1}^n\sum_{(k_1,\ell_1)\in T_{m_1}
\setminus T_{m_1-1}}\sum_{m_2=1}^n\sum_{(k_2,\ell_2)\in T_{m_2}
\setminus T_{m_2-1}}G_{m_1,m_2,k_1,\ell_1,k_2,\ell_2}   \label{eq:eq4.6}\\
&=\sum_{m_1=1}^n\sum_{m_2=1}^n \bigg (\sum_{k_1=-m_1+1}^{m_1}
\sum_{k_2=-m_2+1}^{m_2}\!\!\!
G_{m_1,m_2,k_1,m_1,k_2,m_2} +\!\!\!\sum_{k_1=-m_1+1}^{m_1}\sum_{\ell_2=-m_2+1}^{m_2-1}
\!\!\!G_{m_1,m_2,k_1,m_1,m_2,\ell_2} \nonumber \\
&\phantom{===}+\!\!\!
\sum_{\ell_1=-m_1+1}^{m_1-1}\sum_{k_2=-m_2+1}^{m_2}\!\!\!\!\!\!
G_{m_1,m_2,m_1,\ell_1,k_2,m_2}+\!\!\!\!\sum_{\ell_1=-m_1+1}^{m_1-1}
\sum_{\ell_2=-m_2+1}^{m_2-1}\!\!\!\!\!\!\!
G_{m_1,m_2,m_1,\ell_1,m_2,\ell_2}\bigg ),  \nonumber
\end{align}
where
\begin{align*}
&G_{m_1,m_2,k_1,\ell_1,k_2,\ell_2} \\ 
&\phantom{==}:=\Cov \big((\delta _{m_1,2,k_1-1,\ell_1} 
-\delta _{m_1,2,k_1,\ell_1-1}) \delta _{m_1,2,k_1,\ell_1-1},
(\delta _{m_2,2,k_2-1,\ell_2} -\delta _{m_2,2,k_2,\ell_2-1})
\delta _{m_2,2,k_2,\ell_2-1} \big). 
\end{align*}
By representation \eqref{MArep} of \ $X_{k,\ell}$ \ and definition
of \ $\delta_{m,2,k,\ell}$ \ we have
\begin{alignat*}{2}
\delta_{m,2,k-1,m}&=X_{k-1,m}-
\sum _{i=-m+2}^{k-1}\alpha ^{k-1-i}\varepsilon _{i,m},  &&{-m+2\leq k\leq m},\\
\delta_{m,2,k,m-1}&=X_{k,m-1},  &&{-m+1\leq k\leq m},\\
\delta_{m,2,m,\ell-1}&=X_{m,\ell-1}-
\sum _{j=-m+2}^{\ell-1}(1-\alpha )^{\ell-1-j}\varepsilon _{m,j}, \qquad &&{-m+2\leq 
\ell\leq m-1},\\
\delta_{m,2,m-1,\ell}&=X_{m-1,\ell},  &&{-m+1\leq \ell \leq m-1}. 
\end{alignat*}
Hence, e.g.
\begin{align*}
&\sum_{k_1=-m_1+1}^{m_1}
\sum_{k_2=-m_2+1}^{m_2}
G_{m_1,m_2,k_1,m_1,k_2,m_2} \\
&\phantom{=}=\sum_{k_1=-m_1+1}^{m_1}
\sum_{k_2=-m_2+1}^{m_2}\Cov \bigg( \Big (X_{k_1-1,m_1}-X_{k_1,m_1-1}-
\sum _{i_1=-m_1+2}^{k_1-1}\alpha ^{k_1-1-i_1}\varepsilon _{i_1,m_1}
\Big) X_{k_1,m_1-1}, \\
&\phantom{=\qquad \qquad \qquad \qquad\qquad \qquad}
\Big (X_{k_2-1,m_2}-X_{k_2,m_2-1}-
\sum _{i_2=-m_2+2}^{k_2-1}\alpha ^{k_2-1-i_2}\varepsilon _{i_2,m_2}
\Big) X_{k_2,m_2-1}\bigg) \\
&\phantom{=}=\sum_{k_1=-m_1+1}^{m_1}
\sum_{k_2=-m_2+1}^{m_2} G^{(1)}_{k_1,m_1,k_2,m_2}-
 G^{(2)}_{k_1,m_1,k_2,m_2}-  G^{(2)}_{k_2,m_2,k_1,m_1}+  
G^{(3)}_{k_1,m_1,k_2,m_2},  
\end{align*}
where
\begin{align*}
G^{(1)}_{k_1,m_1,k_2,m_2}&:=\Cov \Big( (X_{k_1-1,m_1}-X_{k_1,m_1-1}
)X_{k_1,m_1-1},(X_{k_2-1,m_2}-X_{k_2,m_2-1})X_{k_2,m_2-1}\Big),\\
G^{(2)}_{k_1,m_1,k_2,m_2}&:=\Cov \bigg( (X_{k_1-1,m_1}-X_{k_1,m_1-1}
)X_{k_1,m_1-1},
X_{k_2,m_2-1}\sum _{i=-m_2+2}^{k_2-1}\alpha
^{k_2-1-i}\varepsilon _{i,m_2} \bigg), \\ 
G^{(3)}_{k_1,m_1,k_2,m_2}&:=\Cov \bigg( X_{k_1,m_1-1}\sum _{i_1=-m_1+2}^{k_1-1}
\alpha ^{k_1-1-i_1}\varepsilon _{i_1,m_1} ,
X_{k_2,m_2-1}\sum _{i_2=-m_2+2}^{k_2-1}\alpha ^{k_2-1-i_2}\varepsilon _{i_2,m_2} \bigg). 
\end{align*}
Thus, Lemma \ref{mmmm}, representation \eqref{MArep} and independence
of the error terms \ $\varepsilon_{i,j}$ \ imply
\begin{align*}
G^{(2)}_{k_1,m_1,k_2,m_2}=&\,\Cov \bigg( X_{k_1-1,m_1}-X_{k_1,m_1-1},
X_{k_2,m_2-1}\bigg)\Cov\bigg( 
X_{k_1,m_1-1},\sum _{i=-m_2+2}^{k_2-1}\alpha
^{k_2-1-i}\varepsilon _{i,m_2} \bigg)\\
&+\Cov \bigg( X_{k_1,m_1-1},X_{k_2,m_2-1}\bigg)\Cov\bigg(
X_{k_1-1,m_1}-X_{k_1,m_1-1},\sum _{i=-m_2+2}^{k_2-1}\alpha
^{k_2-1-i}\varepsilon _{i,m_2} \bigg),\\
G^{(3)}_{k_1,m_1,k_2,m_2}=&\,\Cov \bigg(
X_{k_1,m_1-1},X_{k_2,m_2-1}\bigg) \Cov\bigg(\sum _{i_1=-m_1+2}^{k_1-1}\!\!\!\!
\alpha ^{k_1-1-i_1}\varepsilon _{i_1,m_1} ,\sum _{i_2=-m_2+2}^{k_2-1}\!\!\!\!\alpha
^{k_2-1-i_2}\varepsilon _{i_2,m_2} \bigg)\\
&+\Cov\bigg(X_{k_1,m_1-1},\!\sum _{i_2=-m_2+2}^{k_2-1}\!\!\!\!\alpha
^{k_2-1-i_2}\varepsilon _{i_2,m_2} \bigg)\Cov\bigg(X_{k_2,m_2-1},\!\sum
_{i_1=-m_1+2}^{k_1-1}\!\!\!\!\alpha^{k_1-1-i_1}\varepsilon _{i_1,m_1} \bigg).
\end{align*}
Moreover, using again the independence of the error terms  \
$\varepsilon_{i,j}$ \ one can
easily see that \ $G^{(2)}_{k_1,m_1,k_2,m_2}=0$ \ if \ $m_2>m_1$ \ and
\ $G^{(3)}_{k_1,m_1,k_2,m_2}=0$ \ if \ $m_2\ne m_1$. 
In this way
\begin{align}
  \label{eq:eq4.7}
\sum_{m_1=1}^n\sum_{m_2=1}^n &\sum_{k_1=-m_1+1}^{m_1}
\sum_{k_2=-m_2+1}^{m_2} G^{(3)}_{k_1,m_1,k_2,m_2} \\ =&\sum_{m=1}^n
\sum_{k_1=-m+1}^{m}  
\sum_{k_2=-m+1}^{m} \Cov \big( X_{k_1,m-1},X_{k_2,m-1}\big)\alpha
^{|k_1-k_2|} \sum_{i=0}^{m+k_1\land k_2-3} \alpha^{2 i} \nonumber \\
\leq &\,\frac {C_{\alpha}}{1-\alpha^2}\sum_{m=1}^n \sum_{k_1=0}^{2m-1}
\sum_{k_2=0}^{2m-1} (k_1+k_2)^{1/2}\leq \frac
{3C_{\alpha}}{1-\alpha^2}(n+1)^{7/2}, \nonumber
\end{align}
where the first inequality is a consequence of Lemma
\ref{covbound} and the empty sum is defined to be zero. 

Further, let
\begin{align*}
B^{(1)}_{k_1,m_1,k_2,m_2}:=&\,\Cov\bigg( 
X_{k_1,m_1-1},\sum _{i=-m_2+2}^{k_2-1}\alpha
^{k_2-1-i}\varepsilon _{i,m_2} \bigg), \\
B^{(2)}_{k_1,m_1,k_2,m_2}:=&\,\Cov\bigg(X_{k_1-1,m_1}- 
X_{k_1,m_1-1},\sum _{i=-m_2+2}^{k_2-1}\alpha
^{k_2-1-i}\varepsilon _{i,m_2} \bigg).
\end{align*}
Assuming  \ $m_2< m_1$, \ with the help of representation
\eqref{MArepprob} we obtain 
\begin{align*}
B^{(1)}_{k_1,m_1,k_2,m_2}
=&\!\!\sum _{i=-m_2+2}^{k_1\land k_2-1} \!\!\! {\mathsf P} \big( 
S^{(\alpha)}_{k_1+m_1-m_2-1-i}=k_1-i\big)\alpha
^{k_2-1-i}+\alpha^{k_2-k_1-1}(1-\alpha )^{m_1-m_2-1}  \bone_{\{k_1\leq
  k_2-1\}} , \\
B^{(2)}_{k_1,m_1,k_2,m_2}
=&\!\!\sum _{i=-m_2+2}^{k_1\land k_2-1} \!\! \Big({\mathsf P} \big( 
S^{(\alpha)}_{k_1+m_1-m_2-1-i}=k_1-1-i\big)-{\mathsf P} \big( 
S^{(\alpha)}_{k_1+m_1-m_2-1-i}=k_1-i\big)\Big)\alpha ^{k_2-1-i} \\
& -\alpha^{k_2-k_1-1}(1-\alpha
)^{m_1-m_2-1}  \bone_{\{k_1\leq k_2-1\}} 
\end{align*}
for \ $k_1+m_1\geq 3$, \ otherwise the above quantities are equal to zero.
Hence, using \eqref{eq:eq2.2} one can easily show that for
\ $k_1\leq k_2-1$ \
\begin{align*}
\big|B^{(1)}_{k_1,m_1,k_2,m_2}\big|\leq &\, \alpha^{k_2-k_1-1}(1-\alpha
)^{m_1-m_2-1} + \alpha^{k_2-k_1} \sum _{i=-m_2+2}^{k_1-1}
\frac{D_{\alpha,\alpha}}{(k_1+m_1-m_2-1-i)^{1/2}} \\
\leq &\,H_{\alpha}  \alpha^{k_2-k_1}(k_1+m_1)^{1/2},\\
\big|B^{(2)}_{k_1,m_1,k_2,m_2}\big|\leq &\, \alpha^{k_2-k_1-1}(1-\alpha
)^{m_1-m_2-1} + \alpha^{k_2-k_1} \sum _{i=-m_2+2}^{k_1-1}
\frac{D_{\alpha,\alpha}}{k_1+m_1-m_2-1-i} \\
\leq &\,H_{\alpha}  \alpha^{k_2-k_1}\log (k_1+m_1)
\end{align*}
with some constant \ $H_{\alpha}>0$, \ while for  \ $k_1>
k_2-1$ \  we have
\begin{align*}
\big|B^{(1)}_{k_1,m_1,k_2,m_2}\big|\leq &\,
\frac{D_{\alpha,\alpha}}{(k_1-k_2+m_1-m_2)^{1/2}} \sum
_{i=-m_2+2}^{k_2-1}\alpha ^{k_2-1-i} \leq
\frac{H_{\alpha}}{(k_1-k_2+m_1-m_2)^{1/2}},\\ 
\big|B^{(2)}_{k_1,m_1,k_2,m_2}\big|\leq &\,
\frac{D_{\alpha,\alpha}}{k_1-k_2+m_1-m_2} \sum
_{i=-m_2+2}^{k_2-1}\alpha ^{k_2-1-i} \leq
\frac{H_{\alpha}}{k_1-k_2+m_1-m_2}.
\end{align*} 
Obviously, if \ $m_1=m_2$ \ then
\begin{equation*}
B^{(1)}_{k_1,m_1,k_2,m_2}=0 \quad \text{and} \quad B^{(2)}_{k_1,m_1,k_2,m_2}=\sum
_{i=-m_1+2}^{k_1\land k_2-1} {\mathsf P}
\big(S^{(\alpha)}_{k_1-1-i}=k_1-1-i\big)\alpha^{k_2-1-i}\leq 
\frac{\alpha^{|k_1-k_2|}}{1-\alpha ^2}.
\end{equation*}
In this way, by Lemma \ref{covbound} and Proposition \ref{covdiff}, 
\begin{align}
\sum_{m_1=1}^n&\,\sum_{m_2=1}^n \sum_{k_1=-m_1+1}^{m_1}
\sum_{k_2=-m_2+1}^{m_2}\!\!\big| G^{(2)}_{k_1,m_1,k_2,m_2}\big |\leq \frac
{C_{\alpha}}{1\!-\!\alpha ^2} \sum_{m=1}^n \sum_{k_1=-m+3}^{m}
\sum_{k_2=-m+3}^{m}\!\!\!(k_1\!+\!k_2\!+\!2m)^{1/2}\alpha
^{|k_1-k_2|} \nonumber \\
&+ H_{\alpha} \sum_{m_1=2}^n\sum_{m_2=1}^{m_1-1}\sum_{k_2=-m_2+1}^{m_2}
\sum_{k_1=-m_1+3}^{k_2-1}\!\!\! \alpha^{k_2-k_1}(k_2\!+\!m_1)^{1/2}
\Big(K_{\alpha} +2C_{\alpha}\log(k_2\!+\!m_1)\Big) \label{eq:eq4.8}\\
&+ H_{\alpha}\! \sum_{m_1=2}^n\sum_{m_2=1}^{m_1-1}\sum_{k_1=-m_1+3}^{m_1}
\sum_{k_2=-m_2+1}^{k_1}\!\! \bigg(\frac
{K_\alpha}{(k_1\!-\!k_2\!+\!m_1\!-\!m_2)^{1/2}}+ \frac
{C_{\alpha}(k_1\!+\!k_2\!+\!m_1\!+\!m_2)^{1/2}
}{k_1\!-\!k_2\!+\!m_1\!-\!m_2} \bigg) \nonumber \\ 
\leq &\,\frac {8C_{\alpha}}{(1-\alpha)(1\!-\!\alpha ^2)}\sum_{m=1}^n
m^{3/2} +\frac{H_{\alpha}}{1-\alpha}\sum_{m=2}^n\sum_{k=-m+1}^{m}
m(k+m)^{1/2}\Big(K_{\alpha} +2C_{\alpha}\log(k+m)\Big) \nonumber\\
&+2 H_{\alpha}\sum_{m=2}^n\sum_{k=-m+3}^{m} m(k+m)^{1/2}
\Big(K_{\alpha} +2C_{\alpha}\log(k+m)\Big) \leq
Q_{\alpha}(n+1)^{7/2}\log(n+1) \nonumber
\end{align}
with some constant \ $Q_{\alpha}>0$. \ Inequalities \eqref{eq:eq4.7}
and \eqref{eq:eq4.8} imply
\begin{align*}
\sum_{m_1=1}^n\sum_{m_2=1}^n \sum_{k_1=-m_1+1}^{m_1}&
\sum_{k_2=-m_2+1}^{m_2} G_{m_1,m_2,k_1,m_1,k_2,m_2} \\
&=\sum_{m_1=1}^n\sum_{m_2=1}^n \sum_{k_1=-m_1+1}^{m_1} 
\sum_{k_2=-m_2+1}^{m_2} G^{(1)}_{k_1,m_1,k_2,m_2}+{\mathcal
  O}\big(n^{7/2}\log (n)\big), \nonumber
\end{align*}
and the same can be proved for the remaining three terms of
\eqref{eq:eq4.6}. Hence
\begin{align}
   \label{eq:eq4.9}
\Var \bigg(\sum_{m=1}^n U_m\bigg) = \sum_{m_1=1}^n\sum_{(k_1,\ell_1)\in T_{m_1}
\setminus T_{m_1-1}}\sum_{m_2=1}^n &\sum_{(k_2,\ell_2)\in T_{m_2}
\setminus T_{m_2-1}} G^{(1)}_{k_1,\ell_1,k_2,\ell_2} \\
=&\sum_{(k_1,\ell_1)\in T_n}\sum_{(k_2,\ell_2)\in T_n}
G^{(1)}_{k_1,\ell_1,k_2,\ell_2} +{\mathcal R}_n,    \nonumber
\end{align}
and \ ${\mathcal R}_n={\mathcal O}\big(n^{7/2}\log (n)\big)$. \
Further, Lemma \eqref{mmmm} implies 
\begin{equation*}
G^{(1)}_{k_1,\ell_1,k_2,\ell_2}=
{\mathcal L}^{(1)}_{k_1,\ell_1,k_2,\ell_2}+ {\mathcal
  L}^{(2)}_{k_1,\ell_1,k_2,\ell_2}, 
\end{equation*}
where
\begin{align*}
{\mathcal L}^{(1)}_{k_1,m_1,k_2,m_2}:=&\Cov \big(X_{k_1-1,m_1}-X_{k_1,m_1-1}
,X_{k_2-1,m_2}-X_{k_2,m_2-1}\big)\Cov \big(X_{k_1,m_1-1},X_{k_2,m_2-1}\big)\\
&+\Cov \big( X_{k_1-1,m_1}-X_{k_1,m_1-1},X_{k_2,m_2-1}\big)\Cov
\big(X_{k_2-1,m_2}-X_{k_2,m_2-1},X_{k_1,m_1-1}\big),
\end{align*}
and using the same ideas as in the proof of \eqref{eq:eq3.4} on can
show
\begin{equation}
  \label{eq:eq4.10}
\lim_{n\to\infty} n^{-9/2}\sum_{(k_1,\ell_1)\in  T_n}\sum_{(k_2,\ell_2)\in T_n}
{\mathcal L}^{(2)}_{k_1,\ell_1,k_2,\ell_2} =0.
\end{equation}
Finally, 
\begin{align}
n^{-9/2}&\sum_{(k_1,\ell_1)\in  T_n}\sum_{(k_2,\ell_2)\in T_n}
{\mathcal L}^{(1)}_{k_1,\ell_1,k_2,\ell_2} \nonumber \\ =&
\iint\limits_T\!\!\iint\limits_T \bigg(\sqrt{n}\Cov
\big(Z_{0,1}^{(n)}(s_1,t_1)-Z_{1,0}^{(n)}(s_1,t_1) 
,Z_{0,1}^{(n)}(s_2,t_2)-Z_{1,0}^{(n)}(s_2,t_2)\big)   \label{eq:eq4.11}\\
&\phantom{=====} \times\Cov
\big(Z_{1,0}^{(n)}(s_1,t_1), Z_{1,0}^{(n)}(s_2,t_2)\big)
+\Cov\big(Z_{0,1}^{(n)}(s_1,t_1)-Z_{1,0}^{(n)}(s_1,t_1),Z_{1,0}^{(n)}(s_2,t_2)\big)
\nonumber \\
&\phantom{=====} \times\sqrt{n}\Cov
\big(Z_{0,1}^{(n)}(s_2,t_2)-Z_{1,0}^{(n)}(s_2,t_2),Z_{1,0}^{(n)}(s_1,t_1)\big)\bigg)
{\mathrm d}s_1\,{\mathrm d}t_1\,{\mathrm d}s_2\,{\mathrm d}t_2. \nonumber
\end{align}
With the help of Lemma \ref{covbound} and Proposition \ref{covdiff}
one can easily show that the integrand on the right hand side of
\eqref{eq:eq4.11} can be dominated by \ $K_{\alpha} \big(C_{\alpha}
\sqrt{s_1+t_1+s_2+t_2+1}+K_{\alpha}\big),$ \ which has a finite integral on
\ $T\times T$.  \ Hence, by Proposition \ref{covlim} and dominated
convergence theorem 
\begin{equation*}
\lim_{n\to\infty} n^{-9/2}\sum_{(k_1,\ell_1)\in  T_n}\sum_{(k_2,\ell_2)\in T_n}
{\mathcal L}^{(1)}_{k_1,\ell_1,k_2,\ell_2} =0,
\end{equation*}
which together with \eqref{eq:eq4.9} and \eqref{eq:eq4.10} implies
\eqref{eq:eq4.5}. \proofend

\bigskip
\noindent
{\bf Proof of Proposition \ref{LINDDn}.} \ We have
\begin{equation*}
 \bone_{\left\{\Vert \widetilde
    d_{T_m}-\widetilde d_{T_{m-1}}\Vert\geq \delta \right\}}\leq
\delta^{-2} \Vert \widetilde  d_{T_m}-\widetilde d_{T_{m-1}}\Vert^2,
\end{equation*}
hence to prove the proposition it suffices to show 
\begin{equation*}
\sum_{m=1}^n{\mathsf E}\Big (\big\Vert\widetilde d_{T_m} -\widetilde
d_{T_{m-1}}\big\Vert ^4  \,\Big |\, 
{\mathcal F}_{m-1}\Big) \stoch 0 \qquad \qquad \text{as \ $n\to\infty$},
\end{equation*}
which is a direct consequence of 
\begin{equation*}
n^{-4}\sum_{m=1}^n{\mathsf E}\Big (\big|d_m^{(1)} -d_{m-1}^{(1)}\big| ^4  \,\Big |\, 
{\mathcal F}_{m-1}\Big)\! \stoch\! 0 \quad \text{and} \quad
n^{-5}\sum_{m=1}^n{\mathsf E}\Big (\big|d_m^{(2)} -d_{m-1}^{(2)}\big| ^4  \,\Big |\, 
{\mathcal F}_{m-1}\Big)\! \stoch \!0.
\end{equation*}
However, these statements
have already been proved, see \citet[Section 6, pp. 47-48] {bpz3} and
\citet[Section 4, pp. 31-32] {bpz3}, respectively. \proofend

\bigskip
\noindent
{\bf Acknowledgments.} \  \ Research has been supported by 
the Hungarian  Scientific Research Fund under Grants No. OTKA
T079128/2009 and OTKA NK101680. G. Pap has been partially supported by
the T\'AMOP-4.2.2.A-11/1/KONV project.

\end{document}